\documentclass[reqno]{amsart}
\usepackage[T1]{fontenc}
\usepackage[margin=1in]{geometry}
\usepackage[fontsize=12pt]{fontsize}
\usepackage{amsmath,amssymb,amsthm,mathtools}
\usepackage{needspace}
\usepackage[colorlinks=true,citecolor=blue,linkcolor=blue,urlcolor=blue]{hyperref}
\usepackage[nameinlink,capitalise]{cleveref}

\newtheorem{theorem}{Theorem}[section]
\newtheorem{lemma}[theorem]{Lemma}
\newtheorem{proposition}[theorem]{Proposition}
\theoremstyle{remark}
\newtheorem{remark}[theorem]{Remark}

\DeclareMathOperator{\BiRic}{BiRic}
\DeclareMathOperator{\Ric}{Ric}
\DeclareMathOperator{\Scal}{Scal}
\DeclareMathOperator{\Hess}{Hess}
\DeclareMathOperator{\dist}{dist}
\DeclareMathOperator{\dive}{div}
\DeclareMathOperator{\trace}{tr}
\DeclareMathOperator{\conjrad}{conj}
\DeclareMathOperator{\area}{{\mathcal{H}}^{n-1}}
\newcommand{\bR}{\mathbb R}
\newcommand{\bS}{\mathbb S}

\newcommand{\cF}{\mathcal F}
\newcommand{\cA}{\mathcal A}
\newcommand{\dd}{\,\mathrm d}
\newcommand{\eps}{\varepsilon}
\newcommand{\conjinf}{\conjrad_{\infty}}
\everymath{\displaystyle}
\allowdisplaybreaks[2]

\title[Bi-Ricci curvature and conjugate radius]{A sharp conjugate radius bound for open four-manifolds with positive bi-Ricci curvature}
\author{Ronggang Li}
\address{School of Mathematics and Statistics, Nanjing University of Information Science and Technology, Nanjing, China}
\email{003705@nuist.edu.cn}
\date{}
\hypersetup{pdftitle={A sharp conjugate radius bound for open four-manifolds with positive bi-Ricci curvature},pdfauthor={Ronggang Li}}

\begin{document}
\begin{abstract}
We prove that every complete noncompact Riemannian four-manifold without boundary with bi-Ricci curvature at least $b>0$ has conjugate radius at most $\pi\sqrt{2/b}$. Equality holds on the round cylinder $\bS^3(\sqrt{2/b})\times\bR$. For $b=2$, the bound and the equality model agree with Gromov's prediction for open four-manifolds under a scalar curvature lower bound of $6$.
\end{abstract}
\maketitle

\section{Introduction}

Throughout the paper, an open Riemannian manifold is understood to be connected, complete, noncompact, and without boundary. For $p\in M$, let $\conjrad_g(p)$ denote the infimum of the first conjugate times along unit-speed geodesics issuing from $p$, with value $+\infty$ if there are no conjugate points. The conjugate radius of $(M,g)$ is
\begin{equation}\label{eq:conj-definition}
 \conjrad(M,g)=\inf_{p\in M}\conjrad_g(p).
\end{equation}
Thus $a<\conjrad_g(p)$ guarantees that $\exp_p$ has nonsingular differential on $B_a(0)\subset T_pM$; injectivity is not required.

Zhongmin Shen \cite{Shen01} conjectured that every open manifold with $\Scal_g\ge n(n-1)$ satisfies $\conjrad(M,g)<\pi$. Ge and Zhang \cite{GZ26} proved that if, in addition, $\lambda_0(M)<n$, then $\conjrad(M,g)\le\pi/\sqrt{1-\lambda_0(M)/n}$, where $\lambda_0(M)$ is the bottom of the spectrum of $-\Delta_g$. In particular, their estimate gives $\conjrad(M,g)\le\pi$ when $\lambda_0(M)=0$.

Gromov proposed a dimension-dependent sharp constant. For an open manifold with a positive scalar curvature lower bound, the $k=0$ case of his conjecture \cite{Gro25} predicts
\[
 \bigl(\inf_M\Scal_g\bigr)\conjrad(M,g)^2\le\pi^2(n-1)(n-2).
\]
In dimension four this predicts $\conjrad(M,g)\le\pi$ under $\Scal_g\ge6$. 

Ge and Li \cite{GL26} proved that every open three-manifold with $\Scal_g\ge6$ satisfies $\conjrad(M,g)\le2\pi/3$. To our knowledge, this was the first general bound strictly below $\pi$ in dimension three under the scalar curvature lower bound alone, and it resolves Shen's conjecture in that dimension. Our main theorem gives a sharp bi-Ricci analogue of Gromov's four-dimensional prediction.

For orthonormal tangent vectors $v,w$, define
\begin{equation}\label{eq:biric-definition}
 \BiRic_g(v,w)=\Ric_g(v,v)+\Ric_g(w,w)-\sec_g(v\wedge w).
\end{equation}
We write $\BiRic_g\ge b$ if $\BiRic_g(v,w)\ge b$ for every orthonormal pair. In dimension three, $\BiRic=\Scal/2$. In contrast to a uniform positive Ricci lower bound, a uniform positive bi-Ricci lower bound is compatible with noncompactness, as round cylinders show.

\begin{theorem}\label{thm:main}
Let $(M^4,g)$ be an open Riemannian manifold. If $\BiRic_g\ge b>0$, then
\begin{equation}\label{eq:sharp-main}
 \conjrad(M,g)\le\pi\sqrt{\frac{2}{b}}.
\end{equation}
The constant is optimal: equality holds on $\bS^3(\sqrt{2/b})\times\bR$ with the product metric.
\end{theorem}

The same argument gives $\conjrad(M,g)\le 2\pi/\sqrt{(6-n)b}$ for $3\le n\le5$ when $\BiRic_g\ge b>0$ outside a compact set; see \cref{thm:infinity}. The bound is sharp in dimension four.

In particular, $\BiRic_g\ge2$ gives $\conjrad(M,g)\le\pi$. On the unit cylinder $\bS^3(1)\times\bR$, one has $\min\BiRic_g=2$ and $\Scal_g=6$, so this normalization gives the same sharp model and the same bound as Gromov's prediction. More generally, on $\bS^{n-1}(r_0)\times\bR$, $n\ge3$, one has
\begin{equation}\label{eq:cylinder-values}
 \min\BiRic=\frac{n-2}{r_0^2},\qquad
 \conjrad_g(x)=\pi r_0\quad\text{for every }x.
\end{equation}
We verify these formulas at the end of the paper.

The proof extends the $\mu$-bubble and exponential-map argument of Ge and Li. Keeping the mean curvature terms in the Gauss and Jacobi equations gives spectral Ricci control on the bubbles. The local inradius estimate in \cref{lem:inradius} uses the conformal geodesic argument of Schoen and Yau \cite{SchoenYau83}, in the form developed by Shen and Ye \cite{SY97} in their work on positive bi-Ricci curvature \cite{SY96,SY97}. Antonelli and Xu \cite{AX24} give related spectral diameter and volume comparisons. The inradius estimate forces suitable components of the inverse images under the exponential map to be compact. The degree of these finite coverings then gives the homological obstruction without requiring the bubbles to be simply connected.

\paragraph{\bf AI assistance declaration:}
AI was used during the preparation of this manuscript to assist with English-language writing and editing.

\paragraph{\textbf{Acknowledgments.}}
The author is grateful to Jian Ge and Chuanhuan Li for useful discussions. The author is supported by the Open Project Program of the Key Laboratory of Mathematics and Complex System at Beijing Normal University (Grant No.~202501). 

\section{A local spectral estimate and a homological obstruction}

We use $\Delta=\dive\nabla$. For a two-sided hypersurface with chosen unit normal $\nu$, the second fundamental form and mean curvature are
$A(X,Y)=\langle\nabla_X\nu,Y\rangle$ and $H=\trace A$.

Firstly, we use the following local form of Shen and Ye’s comparison \cite{SY97}, based on the argument of Schoen and Yau \cite{SchoenYau83}

\begin{lemma}\label{lem:inradius}
Let $(N^k,h)$ be a smooth Riemannian manifold, $2\le k\le4$, and let $\Omega\subset N$ be a compact connected smooth domain with nonempty boundary. Suppose that, for a constant $c>0$, a smooth function $u>0$ on a neighborhood of $\Omega$ satisfies
\begin{equation}\label{eq:spectral}
 \Ric_h\ge (\frac{\Delta_hu}{u}+c)\,h.
\end{equation}
Then
\begin{equation}\label{eq:inradius}
 d_\Omega(x,\partial\Omega) \le D(c):=\frac{2\pi}{\sqrt{(5-k)c}}, \qquad \forall x\in\operatorname{Int}\Omega
\end{equation}
where $d_\Omega$ is the intrinsic length distance in $\Omega$. The conclusion is preserved under local isometric pullback of $h$ and $u$.
\end{lemma}

\begin{proof}
Set $f=\log u$ and $\bar h=e^{2f}h$. Compactness gives a shortest $\bar h$-curve in $\Omega$ from $x$ to $\partial\Omega$. Stop it at its first boundary point and parametrize it by $h$-arclength $s\in[0,\ell]$. Its interior lies in $\operatorname{Int}\Omega$ and is a $\bar h$-geodesic with nonnegative fixed-endpoint index form. Write $T$ for its $h$-unit tangent and $f'=\langle\nabla f,T\rangle$. The conformal geodesic equation gives
\begin{equation}\label{eq:conformal-geodesic}
 \nabla_TT=(\nabla f)^\perp,
 \qquad
 \Hess f(T,T)=f''-|\nabla_\perp f|^2.
\end{equation}
Choose an $h$-orthonormal normal frame $E_1,\ldots,E_{k-1}$ parallel for the normal connection along the curve. The conformal connection formula shows that $e^{-f}E_i$ are $\bar h$-parallel and $\bar h$-orthonormal. Summing the index forms of $\varphi E_i$, with $\varphi\in C_c^\infty(0,\ell)$, gives
\begin{align}\label{eq:index-conformal}
  &\sum_{i=1}^{k-1}I_{\bar h}(\varphi E_i,\varphi E_i)  \notag\\
 =&\int_0^\ell e^f\left((k-1)(\varphi'+f'\varphi)^2-\left(\Ric_h(T,T)-(k-2)f''-\Delta_h f\right)\varphi^2\right)\dd s\geq 0.
\end{align}
Indeed, the conformal Ricci formula and \eqref{eq:conformal-geodesic} imply
\begin{align*}
 \Ric_{\bar h}(T,T)
 &=\Ric_h(T,T)-(k-2)\Hess f(T,T)-\Delta_hf
   +(k-2)\bigl((f')^2-|\nabla f|^2\bigr)\\
 &=\Ric_h(T,T)-(k-2)f''-\Delta_hf.
\end{align*}
Substitute $\varphi=e^{-f/2}\psi$ in \eqref{eq:index-conformal} and integrate the $f''\psi^2$ term by parts. Since \eqref{eq:spectral} gives
\[
 \Ric_h(T,T)-\Delta_hf\ge c+|\nabla f|^2,
\]
we obtain
\begin{equation}\label{eq:index-expanded}
\int_0^\ell\left((k-1)(\psi')^2+(3-k)f'\psi\psi'-\frac{5-k}{4}(f')^2\psi^2-|\nabla_\perp f|^2\psi^2-c\psi^2\right)\dd s\geq 0.
\end{equation}
Completing the square yields
\[
(3-k)f'\psi\psi'-\frac{5-k}{4}(f'\psi)^2=-\frac{5-k}{4}\left(f'\psi-\frac{6-2k}{5-k}\psi'\right)^2+\frac{(3-k)^2}{5-k}(\psi')^2.
\]
Using $(k-1)+(3-k)^2/(5-k)=4/(5-k)$, we conclude that
\begin{equation}\label{eq:index-final}
 \int_0^\ell \frac{4}{5-k}(\psi')^2-c\psi^2 \dd s\geq 0.
\end{equation}
By density, this holds for $\psi\in H_0^1(0,\ell)$. Taking $\psi(s)=\sin(\pi s/\ell)$ gives $\ell\le D(c)$, and $d_\Omega(x,\partial\Omega)\le\ell$ proves the claim.
\end{proof}

In particular, when $k=3$ the mixed term in \eqref{eq:index-expanded} vanishes. Thus
\begin{equation}\label{cor:three-dimensional-spectral}
   \int_0^\ell\left(2(\psi')^2-c\psi^2-\frac12(f')^2\psi^2-|\nabla_\perp f|^2\psi^2\right)\dd s \geq 0,
\end{equation}
and the comparison constant is $\pi\sqrt{2/c}$.

We next convert the local estimate into an obstruction in ordinary singular homology with real coefficients.

\begin{proposition}\label{prop:homology}
Let $(M,g)$ be complete, let $Z\subset M$ be closed, and let $\Sigma^k\subset M\setminus Z$ be a closed connected oriented embedded submanifold, $2\le k\le4$. Suppose that a positive smooth function $u$ on $\Sigma$ satisfies
\begin{equation}\label{eq:homology-spectral}
 \Ric_\Sigma\ge\left(\frac{\Delta_\Sigma u}{u}+c\right)g_\Sigma,
 \qquad c>0.
\end{equation}
If $x\in\Sigma$ and
\begin{equation}\label{eq:exp-hypothesis}
 D(c)<a<\conjrad_g(x),\qquad \dist_g(x,Z)>a,
\end{equation}
then the inclusion $\iota:\Sigma\hookrightarrow M\setminus Z$ satisfies
\[
 \iota_*[\Sigma]=0\quad\text{in }H_k(M\setminus Z;\bR).
\]
Consequently, if this class is nonzero and $\dist_g(\Sigma,Z)>D(c)$, then $\conjrad_g(x)\le D(c)$ for every $x\in\Sigma$.
\end{proposition}

\begin{proof}
Since $a<\conjrad_g(x)$, the pullback metric $\widehat g=(\exp_x)^*g$ is smooth and positive definite on $B_a(0)\subset T_xM$, and $\exp_x:B_a(0)\to M$ is a local isometry. The relatively closed embedded submanifold
\[
 \Uparrow_x^{\Sigma}=\exp_x^{-1}(\Sigma)\subset B_a(0)
\]
contains $0_x$. Its induced metric and the function $\widehat u=u\circ\exp_x|_{\Uparrow_x^{\Sigma}}$ satisfy \eqref{eq:homology-spectral}.

Let $r(v)=|v|_{g_x}$. Choose $r_0\in(D(c),a)$ to be a regular value of $r|_{\Uparrow_x^{\Sigma}\setminus\{0_x\}}$, and let $\Omega$ be the connected component of $\Uparrow_x^{\Sigma}\cap\overline B_{r_0}(0)$ containing $0_x$. This is a compact connected smooth manifold, possibly with boundary. If its boundary is nonempty, it lies in $\{r=r_0\}$. Gauss's lemma gives $|\nabla r|_{\widehat g}=1$ away from $0_x$, so every curve in $\Omega$ from $0_x$ to $\partial\Omega$ has length at least $r_0$. Applying \cref{lem:inradius} gives the contradiction
\[
 r_0\le d_\Omega(0_x,\partial\Omega)\le D(c)<r_0.
\]
Thus $\partial\Omega=\varnothing$, and regularity of $r_0$ implies $\Omega\subset B_{r_0}(0)$.

The local diffeomorphism $\pi=\exp_x|_\Omega:\Omega\to\Sigma$ has compact domain. Its image is nonempty, open, and closed, hence equals $\Sigma$. It is therefore a finite covering. With the orientation pulled back from $\Sigma$, its degree is its number of sheets $j\ge1$, and
\begin{equation}\label{eq:degree}
 \pi_*[\Omega]=j[\Sigma].
\end{equation}
The homotopy
\[
 F:[0,1]\times\Omega\longrightarrow M\setminus Z,
 \qquad F(t,v)=\exp_x(tv),
\]
contracts $\iota\circ\pi$ to the constant map at $x$: indeed,
$\dist_g(x,F(t,v))\le t|v|_{g_x}<a$, so its image avoids $Z$. Hence
\[
 0=(\iota\circ\pi)_*[\Omega]=j\,\iota_*[\Sigma].
\]
Real coefficients give the desired vanishing. The final assertion follows by choosing $a$ strictly between $D(c)$ and both $\conjrad_g(x)$ and $\dist_g(x,Z)$.
\end{proof}

\section{Distant stable bubbles}

For a two-sided hypersurface with unit normal $\nu$, we use
$A(X,Y)=\langle\nabla_X\nu,Y\rangle$ and $H=\trace A$.
We first record an algebraic inequality. Let $A$ be a symmetric endomorphism of an $(n-1)$-dimensional Euclidean space, $n\ge3$, and let $e$ be a unit vector. In an orthonormal basis beginning with $e$, write
\[
 A=\begin{pmatrix}a&z^{\mathsf T}\\z&C\end{pmatrix},
\]
where $a=\langle Ae,e\rangle$ and $\trace C=H-a$. Set $C^\circ=C-\frac{H-a}{n-2}\operatorname{Id}$. A direct calculation gives
\[
 |C|^2=|C^\circ|^2+\frac{(H-a)^2}{n-2},
\]
and hence
\begin{align*}
 |A|^2&=a^2+2|z|^2+|C|^2\\
      &=a^2+2|z|^2+|C^\circ|^2+\frac{(H-a)^2}{n-2}\\
      &=|Ae|^2+|z|^2+|C^\circ|^2+\frac{(H-a)^2}{n-2}.
\end{align*}
Therefore,
\begin{align*}
 |A|^2-|Ae|^2+Ha&=|z|^2+|C^\circ|^2+\frac{(H-a)^2}{n-2}+Ha\\
 &=|z|^2+|C^\circ|^2+\frac1{n-2}\left(a+\frac{n-4}{2}H\right)^2+\frac{6-n}{4}H^2.
\end{align*}
In particular,
\begin{equation}\label{eq:Q-bound}
 Q(A,e):=|A|^2-|Ae|^2+H A(e,e) \ge \frac{6-n}{4}H^2.
\end{equation}

\Needspace{12\baselineskip}
\begin{proposition}\label{prop:bubbles}
Let $(M^n,g)$ be an oriented open Riemannian manifold, $3\le n\le5$. Suppose that $\BiRic_g\ge b>0$ outside a compact set $K$. For every $p\in M$, there are arbitrarily large $L$ and closed connected oriented embedded hypersurfaces $\Sigma_L\subset M\setminus K$, carrying positive smooth functions $u_L$, such that
\begin{align}
 \dist_g(p,\Sigma_L)&>L-1,
 & [\Sigma_L]&\ne0\quad\text{in }H_{n-1}(M\setminus\{p\};\bR),
 \label{eq:bubble-position}\\
 \Ric_{\Sigma_L}-\frac{\Delta_{\Sigma_L}u_L}{u_L}g_{\Sigma_L}
 &\ge c_Lg_{\Sigma_L},
 & c_L&=b-\frac{4\pi^2}{(6-n)L^2}>0.
 \label{eq:bubble-spectral}
\end{align}
\end{proposition}

\begin{proof}
Choose a smooth approximation $\rho$ to $\dist_g(p,\cdot)$ such that
\begin{equation}\label{eq:smooth-distance}
 |\rho-\dist_g(p,\cdot)|<1, \qquad |\nabla\rho|\le2.
\end{equation}
see \cite{GW79}. Its sublevel sets are compact. By Sard's theorem, arbitrarily large $L$ can be chosen so that $L,2L,3L$ are regular values of $\rho$. Choose $L$ also large enough that $K\subset\{\rho<L\}$ and $c_L>0$. On the interior of the compact band $N_L=\{L\le\rho\le3L\}$, set
\begin{equation}\label{eq:mu-definition}
 s_L=\frac{\rho-2L}{2},\qquad
 \mu_L=-\frac{4\pi}{(6-n)L}\tan\left(\frac{\pi s_L}{L}\right).
\end{equation}
The function $\mu_L$ tends to $+\infty$ at the inner boundary and to $-\infty$ at the outer boundary. Since $|\nabla s_L|\le1$,
\begin{equation}\label{eq:riccati-profile}
 \frac{6-n}{4}\mu_L^2-|\nabla\mu_L| \ge -\frac{4\pi^2}{(6-n)L^2}.
\end{equation}

For a finite-perimeter set $E$ and an open subset $U\subset M$, the relative measure of the reduced boundary is
\[
 \area(\partial^* E;U)=\sup\left\{\int_E\operatorname{div}_gX\,dV_g:|X|_g\le1,\ \text{$X$ is supported in $U$}\right\},
\]
where the supremum is over smooth compactly supported vector fields. For a smooth boundary, this agrees with its usual area in $U$. We write $\area(\partial^* E)=\area(\partial^* E;N_L^\circ)$.

Let $E_0=N_L\cap\{\rho<2L\}$. Among finite-perimeter sets in
\[
 \cA_L=\{E\subset N_L:E\mathbin{\triangle}E_0\Subset N_L^\circ\},
\]
minimize
\begin{equation}\label{eq:mu-functional}
 \cF_L(E)=\area(\partial^* E)-\int_{N_L}(\chi_E-\chi_{E_0})\mu_L\dd V_g.
\end{equation}
All set equalities are understood up to null sets.
The volume term is a compactly supported difference integral; its two terms are not integrated separately.

For completeness, we describe why a minimizer remains away from the boundary of the band. This is the standard barrier argument for $\mu$-bubbles; see \cite{CL24}. In sufficiently small regular collars of the inner and outer boundaries, respectively, put
\[
 \eta_{\mathrm{in}}=\frac{\nabla\rho}{|\nabla\rho|},\qquad
 \eta_{\mathrm{out}}=-\frac{\nabla\rho}{|\nabla\rho|}.
\]
Their divergences are bounded. Thus one can choose $0<\delta<L$ such that the collars
\[
 C_{\mathrm{in}}=\{L\le\rho<L+\delta\},\qquad
 C_{\mathrm{out}}=\{3L-\delta<\rho\le3L\}
\]
have regular interior interfaces and satisfy
\[
 \mu_L-\dive\eta_{\mathrm{in}}\ge1,\qquad
 \mu_L+\dive\eta_{\mathrm{out}}\le-1.
\]
For either collar $C$, let $B=C\cap\partial N_L$ and let $\eta$ be the corresponding field. The areas of reduced boundaries are relative to $N_L^\circ$, whereas $\area(B)$ denotes the ordinary area of $B$. Suppose that $A$ fills a relative neighborhood of $B$ in $N_L$. The field points outward from the collar on its interior interface and inward to $N_L$ on $B$. Gauss--Green gives
\[
 \int_C\dive\eta\dd V_g=\area(\partial^*C)-\area(B),\qquad
 \int_{A\cap C}\dive\eta\dd V_g\le\area\left(\partial^*(A\cap C)\right)-\area(B).
\]
Subtracting and using submodularity of perimeter, we obtain
\begin{equation}\label{eq:collar-perimeter}
 \area\left(\partial^*(A\cup C)\right)-\area(\partial^*A)\le\area(\partial^*C)-\area\left(\partial^*(A\cap C)\right)\le\int_{C\setminus A}\dive\eta\dd V_g.
\end{equation}
Apply this with $A=E$ at the inner boundary and $A=N_L\setminus E$ at the outer boundary. It follows that
\begin{align}
 \cF_L(E\cup C_{\mathrm{in}})-\cF_L(E)
 &\le\int_{C_{\mathrm{in}}\setminus E}
      (\dive\eta_{\mathrm{in}}-\mu_L)\dd V_g\le0,
 \label{eq:inner-barrier}\\
 \cF_L(E\setminus C_{\mathrm{out}})-\cF_L(E)
 &\le\int_{E\cap C_{\mathrm{out}}}
      (\dive\eta_{\mathrm{out}}+\mu_L)\dd V_g\le0.
 \label{eq:outer-barrier}
\end{align}
Hence the infimum is unchanged if competitors must fill $C_{\mathrm{in}}$ and avoid $C_{\mathrm{out}}$.

The sign of $\mu_L$ relative to $E_0$ gives
\begin{equation}\label{eq:perimeter-control}
 \cF_L(E)=\area(\partial^* E)+\int_{E\triangle E_0}|\mu_L|\dd V_g\ge \area(\partial^* E).
\end{equation}
A minimizing sequence in this fixed-collar class therefore has bounded perimeter, and its difference from $E_0$ is supported in a fixed compact core where $\mu_L$ is bounded. BV compactness, lower semicontinuity of perimeter, and $L^1$ continuity of the volume term yield a minimizer $E_*$. By \eqref{eq:inner-barrier}--\eqref{eq:outer-barrier}, it also minimizes over the original class $\cA_L$, so all compactly supported interior variations are allowed. The regularity theory for prescribed mean curvature minimizers gives a smooth embedded boundary in our dimensions \cite{CL24}. Gluing $E_*$ to $\{\rho<L\}$ produces a compact smooth region $W_L$ with
\begin{equation}\label{eq:W-location}
 \{\rho\le L\}\subset W_L\subset\{\rho<3L\},
 \qquad \partial W_L\subset\{L<\rho<3L\}.
\end{equation}
Neither the band nor $\partial W_L$ needs to be connected.

Let $\Sigma$ be a connected component of $\partial W_L$, with outward unit normal $\nu$. The first and second variations of \eqref{eq:mu-functional} give $H=\mu_L$ and
\begin{equation}\label{eq:mu-stability}
 \int_\Sigma\left[|\nabla_\Sigma\phi|^2-\bigl(\Ric_M(\nu,\nu)+|A|^2+\nu(\mu_L)\bigr)\phi^2\right]\dd A\ge0
\end{equation}
for all $\phi\in C^\infty(\Sigma)$. Thus the first eigenvalue $\lambda$ of
\[
 J=-\Delta_\Sigma-\bigl(\Ric_M(\nu,\nu)+|A|^2+\nu(\mu_L)\bigr)
\]
is nonnegative. A positive first eigenfunction $u$ satisfies
\begin{equation}\label{eq:jacobi-eigenfunction}
 -\frac{\Delta_\Sigma u}{u}
 =\Ric_M(\nu,\nu)+|A|^2+\nu(\mu_L)+\lambda.
\end{equation}
For a unit tangent vector $e$, the Gauss equation reads
\[
 \Ric_\Sigma(e,e)=\Ric_M(e,e)-\sec_M(e\wedge\nu)+HA(e,e)-|Ae|^2.
\]
Combining this with \eqref{eq:jacobi-eigenfunction} gives
\begin{equation}\label{eq:gauss-jacobi}
 \Ric_\Sigma(e,e)-\frac{\Delta_\Sigma u}{u}
 =\BiRic_M(e,\nu)+Q(A,e)+\nu(\mu_L)+\lambda.
\end{equation}
Since the band is disjoint from $K$, equations \eqref{eq:Q-bound} and \eqref{eq:riccati-profile}, together with $H=\mu_L$, imply
\[
 \Ric_\Sigma(e,e)-\frac{\Delta_\Sigma u}{u}
 \ge b+\frac{6-n}{4}\mu_L^2-|\nabla\mu_L|\ge b-\frac{\pi^2}{\frac{6-n}{4}L^2}=c_L.
\]
Thus every component of $\partial W_L$ satisfies \eqref{eq:bubble-spectral}.

For large $L$, a small normal ball $B_\eps(p)$ lies in $\operatorname{Int}W_L$. The oriented boundary of $W_L\setminus B_\eps(p)$ gives
\begin{equation}\label{eq:boundary-class}
 [\partial W_L]=[\partial B_\eps(p)]
 \quad\text{in }H_{n-1}(M\setminus\{p\};\bR).
\end{equation}
Since $M$ is connected, oriented, and noncompact, $H_n(M;\bR)=0$. Excision and the long exact sequence of the pair therefore give an injection
\[
 H_n(M,M\setminus\{p\};\bR)\simeq\bR
 \lhook\joinrel\longrightarrow H_{n-1}(M\setminus\{p\};\bR),
\]
whose image contains $[\partial B_\eps(p)]$ as a nonzero class \cite{Hat02}. By \eqref{eq:boundary-class}, at least one of the finitely many components of $\partial W_L$ has nonzero class. Choose this component as $\Sigma_L$ and its positive eigenfunction as $u_L$. Equations \eqref{eq:smooth-distance} and \eqref{eq:W-location} give $\dist_g(p,\Sigma_L)>L-1$.
\end{proof}

\begin{remark}\label{rem:mu-smallness}
The relevant small error is $\frac{6-n}{4}\mu_L^2-|\nabla\mu_L|\ge-o(1)$ as $L\to\infty$. The function $\mu_L$ itself diverges at the ends of the band. In dimension four, writing $q_4=(6-4)/4$, we have
\[
 q_4=\frac12,\qquad
 \mu_L=-\frac{2\pi}{L}\tan\left(\frac{\pi s_L}{L}\right),\qquad
 c_L=b-\frac{2\pi^2}{L^2}.
\]
The positive $H^2/2$ term in \eqref{eq:Q-bound} absorbs the derivative of this profile with an error tending to zero.
\end{remark}

\section{Conjugate radius bounds and sharpness}

For an open Riemannian manifold, define
\begin{equation}\label{eq:conj-infinity}
 \conjinf(M,g)
 :=\lim_{R\to\infty}\ \inf_{x\in M\setminus B_g(p,R)}\conjrad_g(x)
 =\liminf_{x\to\infty}\conjrad_g(x).
\end{equation}
The limit exists in $[0,+\infty]$, is independent of $p$, and is at least $\conjrad(M,g)$.

\begin{theorem}\label{thm:infinity}
Let $(M^n,g)$ be an open Riemannian manifold, $3\le n\le5$. Suppose that $\BiRic_g\ge b>0$ outside a compact set $K$. For every $\eps>0$, there exists $A=A(n,b,\eps)>0$ such that, for any $p\in M$ and all sufficiently large $R$, there is a compact connected smooth region $\Omega_R$ with
\[
 \overline B_g(p,R)\subset\operatorname{Int}\Omega_R,\qquad \Omega_R\subset B_g(p,R+A),
\]
and
\[
 \sup_{x\in\partial\Omega_R}\conjrad_g(x)\le\frac{2\pi}{\sqrt{(6-n)b}}+\eps.
\]
Consequently, for every continuous proper ray $\gamma:[0,\infty)\to M$,
\[
 \liminf_{t\to\infty}\conjrad_g(\gamma(t))\le\frac{2\pi}{\sqrt{(6-n)b}}.
\]
In particular,
\begin{equation}\label{eq:general-bound}
 \conjinf(M,g)\le \frac{2\pi}{\sqrt{(6-n)b}}.
\end{equation}
The same upper bound holds for $\conjrad(M,g)$, and the constant is optimal when $n=4$.
\end{theorem}

\begin{proof}
Fix $\eps>0$ and choose $L$ sufficiently large that $c_L>0$ and
\[
 D_L=D(c_L)=\frac{2\pi}{\sqrt{(6-n)c_L}}\le\frac{2\pi}{\sqrt{(6-n)b}}+\eps.
\]
Set $A=2L+3+D_L$. Fix $p\in M$ and take $R>D_L$ large enough that $K\subset B_g(p,R)$. Choose $\rho$ as in \eqref{eq:smooth-distance}. By Sard's theorem we may choose $a\in(R+1,R+2)$ such that $a,a+L,a+2L$ are regular values of $\rho$.

Apply the construction of $W_L$ in the proof of \cref{prop:bubbles} with $\rho-a+L$ in place of $\rho$. The shifted function has compact sublevel sets, the same gradient bound, and regular values $L,2L,3L$. The band is now $\{a\le\rho\le a+2L\}$, which lies outside $K$. The estimates \eqref{eq:riccati-profile} and \eqref{eq:bubble-spectral} therefore hold with the same $c_L$. The construction through \eqref{eq:bubble-spectral} uses the Riemannian volume density and the outward unit normals of the boundary, so it also applies when $M$ is nonorientable. By \eqref{eq:W-location} for the shifted function and the distance estimate for $\rho$, we have
\[
 \overline B_g(p,R)\subset\operatorname{Int}W_L,\qquad W_L\subset B_g(p,R+2L+3).
\]
Let $Z=\{p\}\cup(M\setminus B_g(p,R+A))$. Since $R>D_L$ and $A=2L+3+D_L$, the triangle inequality gives $\dist_g(\partial W_L,Z)>D_L$.

Let $S_L$ be the union of the components of $\partial W_L$ on which
\begin{equation}\label{eq:pointwise-on-bubbles}
 \conjrad_g(x)\le D_L
\end{equation}
holds at every point. We show that the component of $M\setminus S_L$ containing $p$ has compact closure in $B_g(p,R+A)$.

When $M$ is oriented, let $\pi:\widetilde M\to M$ be the identity map. Otherwise, take $\pi$ to be the orientable double cover with the complete pullback metric. Write $\widetilde W_L=\pi^{-1}(W_L)$ and $\widetilde Z=\pi^{-1}(Z)$. Each connected component of $\partial\widetilde W_L$ covers a component of $\partial W_L$ and inherits \eqref{eq:bubble-spectral}. Its distance from $\widetilde Z$ is greater than $D_L$, since projection preserves the length of curves. Geodesics and Jacobi fields lift, so the pointwise conjugate radius is preserved by $\pi$.

Let $c$ be a finite fundamental chain of the compact oriented region $\widetilde W_L$, for instance one given by a triangulation. Its boundary is the sum of the fundamental cycles of the components of $\partial\widetilde W_L$ with their boundary orientations. Each component outside $\pi^{-1}(S_L)$ contains a lift of a point where \eqref{eq:pointwise-on-bubbles} fails. By \cref{prop:homology}, its fundamental cycle bounds a finite real chain in $\widetilde M\setminus\widetilde Z$. Subtract all of these filling chains from $c$ and denote the resulting $n$-chain by $c'$. Then $\partial c'$ is supported in $\pi^{-1}(S_L)$ and $\pi(\operatorname{supp}c')$ is a compact subset of $B_g(p,R+A)$.

For $y\in\widetilde M\setminus\pi^{-1}(S_L)$, the chain $c'$ defines a local class in $H_n(\widetilde M,\widetilde M\setminus\{y\};\bR)$. Its coefficient relative to the orientation class is locally constant there. Indeed, on a small closed coordinate ball disjoint from $\pi^{-1}(S_L)$, the classes at interior points come from the same relative class of $c'$ for the ball, by excision \cite{Hat02}. This coefficient is $1$ at each lift of $p$, since the fillings avoid those lifts and $p\in\operatorname{Int}W_L$. It is $0$ outside the support of $c'$. Every path in $M\setminus S_L$ starting at $p$ lifts to a path on which the coefficient stays equal to $1$. Thus the component $U$ of $M\setminus S_L$ containing $p$ is contained in $\pi(\operatorname{supp}c')$.

The closed ball $\overline B_g(p,R)$ is connected and disjoint from $S_L$, so it lies in $U$. Each component of $S_L$ is two-sided, with two connected sides in a collar. A side meeting $U$ is contained in $U$. The components incident to $U$ on exactly one side form the boundary of $\overline U$, and those incident on both sides lie in its interior. Hence $\Omega_R=\overline U$ is a compact connected smooth region with $\partial\Omega_R\subset S_L$, and the claimed bound follows from \eqref{eq:pointwise-on-bubbles}.

For a proper ray, take $p=\gamma(0)$. The ray eventually leaves each compact $\Omega_R$, so it meets $\partial\Omega_R$ at a distance greater than $R$ from $p$. The corresponding parameters tend to infinity as $R\to\infty$. Letting $\eps\downarrow0$ gives the bound along the ray. To obtain \eqref{eq:general-bound}, take points on $\partial\Omega_R$ as $R\to\infty$ and let $\eps\downarrow0$ again. The global bound follows from $\conjrad(M,g)\le\conjinf(M,g)$. The cylinder calculation below proves sharpness in dimension four.
\end{proof}

For a unit-speed minimizing geodesic ray from $p$, $\dist_g(p,\gamma(t))=t$. Hence for each $\eps>0$ and each sufficiently large $R$, the ray meets $\partial\Omega_R$ at some parameter in $[R,R+A]$, where the pointwise conjugate radius satisfies the stated bound.

\begin{proof}[Proof of \cref{thm:main} and of sharpness]
The inequality follows from \cref{thm:infinity} with $n=4$ and an empty exceptional set. To check sharpness, consider $\bS^{n-1}(r)\times\bR$. For an orthonormal pair $v,w$, set
\[
 a=\langle v,\partial_t\rangle,\qquad
 \beta=\langle w,\partial_t\rangle,\qquad \theta=a^2+\beta^2\in[0,1].
\]
The product curvature formulas give
\[
 \Ric(v,v)=\frac{n-2}{r^2}(1-a^2),\qquad
 \Ric(w,w)=\frac{n-2}{r^2}(1-\beta^2),\qquad
 \sec(v\wedge w)=\frac{1-\theta}{r^2}.
\]
Consequently,
\[
 \BiRic(v,w)=\frac{2n-5-(n-3)\theta}{r^2}\ge\frac{n-2}{r^2},
\]
with equality for a two-plane containing $\partial_t$.

A unit-speed product geodesic has constant speed $s\in[0,1]$ in the sphere factor. If $s>0$, its first conjugate time is $\pi r/s$; if $s=0$, it has no conjugate point. Taking the infimum over directions yields $\conjrad_g(x)=\pi r$ at every point. This proves \eqref{eq:cylinder-values}. For $n=4$ and $r=\sqrt{2/b}$, equality holds in \eqref{eq:sharp-main}, and also in \eqref{eq:general-bound}.
\end{proof}

\begin{remark}
The estimate concerns the infimum of first conjugate times over directions. It does not bound the first conjugate time along every geodesic, as the line geodesics in the cylinder show. The equality example proves optimality of the constant; no classification of equality cases is asserted.
\end{remark}

\end{document}